\documentclass[review]{elsarticle}

\usepackage[utf8]{inputenc}
\usepackage{amsmath}
\usepackage{amsfonts}
\usepackage{amssymb}
\usepackage{graphicx}
\usepackage{fullpage}
\usepackage{hyperref}
\usepackage{todonotes}
\usepackage[notcite,color,notref]{showkeys}
\usepackage[titletoc,title]{appendix}
\usepackage[mathscr]{eucal}
\hypersetup{pdfauthor={Name}}

\newtheorem{thm}{Theorem}[section]
\newtheorem{lem}[thm]{Lemma}

\newtheorem{remark}{Remark}[section]
\newtheorem{assum}{Assumption}[section]
\newcommand{\proof}{\noindent {\bf Proof.}~}

\newcommand{\bR}{\mathbb{R}}
\newcommand{\bZ}{\mathbb{Z}}

\newcommand{\bu}{\mathbf{u}}
\newcommand{\bv}{\mathbf{v}}
\newcommand{\bx}{\mathbf{x}}

\newcommand{\cB}{\mathcal{B}}
\newcommand{\cF}{\mathcal{F}}

\newcommand{\cI}{\mathcal{I}}

\newcommand{\cU}{\mathcal{U}}
\newcommand{\cV}{\mathcal{V}}
\newcommand{\cX}{\mathcal{X}}

\begin{document}

\title{Transitions in echo index and dependence on input repetitions}

\author[1]{Peter Ashwin\corref{cor1}}\ead{p.ashwin@exeter.ac.uk}

\author[1,2]{Andrea Ceni}\ead{andrea.ceni.19@gmail.com}

\cortext[cor1]{Corresponding author}
\address[1]{Department of Mathematics, University of Exeter, Exeter EX4 4QF, UK}
\address[2]{Department of Computer Science, University of Pisa, Largo Bruno Pontecorvo, 3 - 56127, IT}

\begin{abstract}
The \emph{echo index} counts the number of simultaneously stable asymptotic responses of a nonautonomous (i.e. input-driven) dynamical system. It generalizes the well-known \emph{echo state property} for recurrent neural networks - this corresponds to the echo index being equal to one. In this paper, we investigate how the echo index depends on parameters that govern typical responses to a finite-state ergodic external input that forces the dynamics.
We consider the echo index for a nonautonomous system that switches between a finite set of maps, where we assume that each map possesses a finite set of hyperbolic equilibrium attractors. 
%%%%%
We find the minimum and maximum repetitions of each map are crucial for the resulting echo index.
Casting our theoretical findings in the RNN computing framework, we obtain that for small amplitude forcing the echo index corresponds to the number of attractors for the input-free system, while for large amplitude forcing, the echo index reduces to one. The intermediate regime is the most interesting; in this region the echo index depends not just on the amplitude of forcing but also on more subtle properties of the input. 
\end{abstract}


\maketitle

\newpage

\tableofcontents

\section{Introduction}

One of the most pressing questions for artificial intelligence systems, is whether one can understand and query the reasons behind a decision made by such a system. The difficulty of answering this {\em explainability problem} \cite{Explain2022} is reflected in the fact that a trained neural network is commonly referred to as a black box. This suggests it is important to try and understand the functioning of neural networks in decision making under input - it is important to {\em open the black box} \cite{sussillo2013opening} and nonlinear dynamics gives tools that can be used for this. As an example, \cite{ceni2020interpreting} show that excitable network attractors can be used to understand function and malfunction of trained RNNs for certain tasks.

Recurrent neural networks (RNN) such as echo state networks \cite{jaeger2001echo,lukovsevivcius2009reservoir} can retain memory of internal states. In this case, it is important to view the system in an input-driven context \cite{manjunath2012theory}. A useful criterion for successful computation is the Echo State Property (ESP) \cite{jaeger2001echo}, which holds if there is asymptotic loss of information about the internal state of the system and only the input is important to determine the output. However, as discussed in \cite{manjunath2013echo}, the presence of the ESP depends not just on system but on the particular input considered. As the input streams to an input-driven RNN will never be fully deterministic, this means there is no guarantee that ESP will be satisfied in practise. Recent work has highlighted that RNNs that do not have ESP may still be a useful model for understanding errors in neural networks  \cite{ceni2020echo}, or in order to design multifunction RNNs that switch between different tasks \cite{ceni2021phd,multifunction2022recurrent}.

In this paper we develop ideas in \cite{ceni2020echo} in more depth by considering parametrizations of input by repetition properties of the inputs. The remainder of this section discusses shift dynamics and introduces the echo index for discrete time dynamical systems. We discuss parametrization of inputs by symbol repetition. Section~\ref{sec:minecho} presents some conditions to give bounds on echo index and the ESP.
Section~\ref{sec:rnns} turns to a numerical example of an echo state RNN where we show how the echo index changes depending on min-max block-length and parameters of the input to the RNN. Section~\ref{sec:conclusions} concludes with a discussion, including barriers to strengthening the theoretic results in Section~\ref{sec:rnns}.

\subsection{Input-driven dynamics and shift dynamics}

We consider properties of a discrete time input-driven dynamical system \cite{manjunath2012theory} of the form
\begin{equation}
\label{eq:driven_dynamics}
x[k+1] = f(x[k],u[k]),
\end{equation}
where time is indexed by $k\in\bZ$, states are $x[k]\in X\subset \bR^n$ and $u[k]\in U\subset \bR^p$ is an input sequence. We denote sequences in bold font, i.e. 
$$
\bx = \{ x[k] \}_{k \in \bZ}, ~~\bu = \{ u[k] \}_{k \in \bZ}, 
$$
and space of sequences in calligraphic font, i.e. $\cX = X^\bZ$, $\cU = U^\bZ$, so that $\bx \in \cX$, $\bu \in \cU$. We assume that $X$ and $U$ are compact sets and assume that $f:X\times U\rightarrow X$ is an update function that is continuous in both arguments, so that a forward orbit is defined. Note that for a given input sequence $u[k]$, (\ref{eq:driven_dynamics}) can be thought of as a nonautonomous dynamical system \cite{kloeden2011nonautonomous}
\begin{equation}
\label{eq:nad_dynamics}
x[k+1] = F(x[k],k),
\end{equation}
where $F:X\times \bR\rightarrow X$ is defined by $F(x,k)=f(x,u[k])$.  The system (\ref{eq:driven_dynamics}) can also be viewed as a skew product dynamical system \cite{kloeden2011nonautonomous} over shift dynamics on the input sequence $\bu\in \cU$. 

We consider the special case where $U$ is taken from a finite set of $M$ input values. Without loss of generality we denote $U=\{0,\ldots,M-1\}$ and call elements of $U$ {\em symbols}. We recall some standard concepts from the symbolic dynamics of shifts; see for example \cite{adler1992symbolic,lind1995introduction,bruin2022topological} for background and more details. We define the shift operator $\sigma$ by
$$
[\sigma({\bf u})][k]=u[k+1].
$$
such that $\sigma:\cU\rightarrow \cU$. Consider a given subset of input sequences $\cV \subseteq \cU$. 
We say $\cV$ is \emph{shift-invariant} if $ \sigma(\cV) = \cV $ and call $(\sigma,\cU)$ the full shift on $U$. We consider the product topology on $\cV$ induced by the metric
\begin{equation}
d(\bu,\bv)=\sum_{k\in\bZ} \frac{d_U(u[k],v[k])}{2^{|k|}}
\label{eq:metric}
\end{equation}
where $d_U(u,v)$ is a metric on $U$.
Note that $\sigma$ acts continuously under such a topology.

If we assume that $\cV$ is shift-invariant and closed then one can lift the nonautonomous system \eqref{eq:driven_dynamics} to a continuous \emph{autonomous} dynamical system on an extended space
$\cF: X \times \cV \rightarrow X \times \cV$. This is given by 
\begin{equation}
\label{eq:extended-autonomous}
    [x[k+1],\bu[k]]=\cF(x[k],{\bf u}[k]) =  \left( f(x,u[0]) ,  \sigma({\bf u}[k])    \right)
\end{equation}
where we write $\bu[k]=\sigma^{k}\bu$. The skew product nature of (\ref{eq:extended-autonomous}) means that composition can be written
\begin{equation}
\label{eq:with-cocycle}
\cF^n(x,{\bf u}) = \left(\Phi_{n,{\bf u}}(x), \sigma^n(\bu)  \right)
\end{equation}
for all $n\in\bZ_0^+$, where $\Phi$ is a cocycle, namely $\Phi_{0,{\bu}}$ is the identity map on $X$ and 
$\Phi_{n+k,\bu}= \Phi_{n,\sigma^k(\bu)}\circ \Phi_{k,\bu}$
for any $n,k\in\bZ^+_0$. We write the forward orbit of $x[0]$ driven by the input sequence $\bu$ in terms of this cocycle, as
\begin{equation}
    x[k]=\Phi_{k,\bu}(x[0]).
\end{equation}
for any $k\in\bZ_0^+$. This can be extended to negative $k$ if $f$ is invertible.

There are many choices of $\cV\subset \cU$ that may be used to characterise a set of possible input sequences, for convenience we only consider closed and shift-invariant subsets. Given any $\bu\in\cU$ one can define a closed invariant subset in terms of its orbit closure
$$
\cU(\bu)=\overline{\{\sigma^{k}(\bu)~:~k\in\bZ\}}
$$
where $\overline{\cV}$ denotes the closure of $\cV$ in the topology of (\ref{eq:metric}).

An important set of closed and shift-invariant subsets are {\em subshifts of finite type} defined as follows: For any $m\geq k$ we write  $u[k,m]=(u[k],u[k+1],\ldots,u[m])$ to denote a finite string or {\em word} of symbols. Given a finite set $\cB$ of words we define
$$
\cU_{\cB}=\{\bu\in\cU~:~u[k,m]\cap \cB=\emptyset \mbox{ for all } k,m \in \bZ \mbox{ with } k\leq m\},
$$
namely the set of sequences that contain no word in the set $\cB$. This is called a \emph{subshift of finite type} with \emph{forbidden words} $\cB$. Without loss of generality we can assume that $\cB$ is \emph{minimal} in the sense there is no proper subset $\cB'\subset \cB$ such that
$\cU_{\cB}=\cU_{\cB'}$.

If $\cU_{\cB}$ has a set of forbidden words of length at most $k+1$ then we say it is a \emph{$k$-step subshift}: knowing $k$ consecutive symbols provides the only constraints on the next symbol.
Note that by considering a higher block shift on overlapping words of length $m$
% $$
% \bu^{m^+}=\cdots\,u[1-m^+,0]\,u[2-m^+,1]\,u[3-m^+,2]\,\cdots
% $$
it is possible to express such a $\cU$ as a 1-step subshift of finite type \cite{bruin2022topological}.

We write system (\ref{eq:driven_dynamics}) in the form of an iterated function system
\begin{equation}
\label{eq:switching_dynamics}
x[k+1] = f_{u[k]}( x[k]),
\end{equation}
for $k\in\bZ$, where \eqref{eq:switching_dynamics} describes the dynamics of switching between one of $M$ of autonomous dynamical systems.

\subsection{Local attractors and the echo index}

We recall a nonautonomous notion of local attractor of \eqref{eq:driven_dynamics} from \cite{ceni2020echo} for fixed input sequence; this was used to proposed a generalization of the notion of ESP for input-driven systems with finitely many local attractors.
For a given input sequence $\bv\in\cU$, we say $\bx=\{x[k]\}_{k\in\bZ}$ is an {\em entire solution} if it is a trajectory for the input $\bv$ in forwards and backwards time; i.e. it satisfies (\ref{eq:driven_dynamics}) for all $k\in\bZ$. 

Note that we do not require $f_i$ to be invertible or surjective, and hence given some $x[0]\in X$ there may be (a) multiple choices for $x[-1]$ or (b) no choice for $x[-1]$ that lies within $X$. Note also that whether (a) or (b) hold will typically depend on $v[k]$ for $k<0$.

Note that there will always be a point $x[0]$ that is on an entire solution $\bx$. In particular, the set
$$
X_{-n, 0}(\bv):=\Phi_{n,\sigma^{-n}(\bv)}(X)
$$
is well-defined, compact and non-empty as $f(X,v)\subset X$ is a continuous image of a compact set. Moreover
$$
X_{-n-1, 0}(\bv)=\Phi_{n,\sigma^{-n}(\bv)}\circ f_{v[-n-1]} (X)\subset X_{-n, 0}(\bv)
$$
and hence the set
$$
X_0(\bv):=\bigcap_{n>0} X_{-n, 0}(\bv)
$$
consists of points that lie on entire solutions.

An entire solution $\bx$ is {\em globally (pullback) attracting} for input $\bv$ if 
$$
\lim_{n\rightarrow \infty} h( X_{-n, 0}(\bv), x[0])=0
$$
where $h(A,B)=\sup_{y\in A} \inf_{z\in B} d(y,z)$.
Following \cite{manjunath2013echo} say (\ref{eq:driven_dynamics}) for a given input $\bv$ has the {\em ESP} if it has a unique entire solution $\bx$ that is globally attracting for this input. An equivalent condition is that there is single point in $X_0(\bv)$.

However, just as an autonomous dynamical system may be multistable, there may be more than one locally attracting entire solution. Moreover, the system (\ref{eq:driven_dynamics}) may have the echo state property for some inputs but not for others. This is explored in \cite{ceni2020echo} where we define a notion of echo index as the smallest number of uniform attraction entire solutions (UAES) that attract almost all initial states of the system.

If there are a number $m$ of UAESs $\{ \bx_1, \ldots, \bx_m \}$ such that they decompose\footnote{apart from a subset of zero Lebesgue measure.} the whole phase space in sets that are uniformly attracted to them, then we say the system has {\em echo index} $m$ and write 
$\cI(\bv)=m$; see \cite[Definitions 3.3 and 3.4]{ceni2020echo} for formal definitions. Note that the echo index is shift-invariant, that is $\cI(\sigma({\bf v}))=\cI({\bf v})$. We note that the echo index $\cI$ may take a number of values on a given closed shift-invariant set $\cV\subset \cU$. In particular, we say $\cV$  has \emph{consistent echo} $n$ if $\cI(\bv)=n$ for all $\bv\in\cV$.

\section{Sufficient conditions to determine minimum echo index}
\label{sec:minecho}

We show that certain assumptions on the behaviour of individual maps and input sequences can be used to guarantee minimum echo index in terms of symbol repetitions for the input.

\subsection{Subshifts with min-max repetitions of input sequences}
\label{sec:inputs}

Suppose we have a finite set of symbols in $U$ and consider a subshift $U_{\cB}$ of finite type defined by a set of prohibited words $\cB$. 
We consider a particular class of subshifts $\cU_{\cB}$ that are characterised by minimum and maximum numbers of a repetition.

We say there is {\em min repetition $m^-_i$ of symbol $i$} if
$\cB$ contains all words of the form
$$
j\,i^{m}\,k
$$
for all $m<m^-_i$ and all $j,k\neq i$. Similarly, we say there is {\em max repetition $m^+_i$ of symbol $i$}
if $\cB$ contains the word
$$
i^{m^+_i+1}.
$$
We say $m_i^+=\infty$ if there is no such word in $\cB$. We say the subshift of finite type $\cU_{\cB}$ has {\em min-max repetitions} $m^{\pm}_i$ in this case. 

\begin{figure}
\centering
\includegraphics[width=8cm]{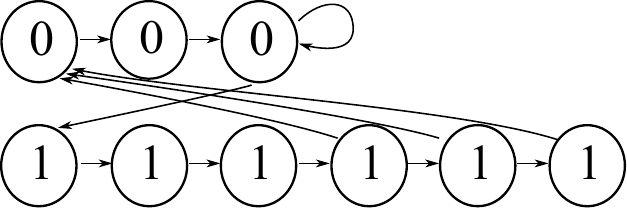}
\caption{Graphical representation of a subshift that has minimum $3$ repeats of $0$ and between $4$ and $6$ repeats of $1$. For this case we have min-max repetitions $m_0^-=3$, $m_0^+=\infty$, $m_1^-=4$ and $m_1^+=6$.}
\label{fig:minmaxexample}
\end{figure}

For example, for the two symbols $0$ and $1$ consider the subshift with prohibited words
$$
\cB=\{010,0110,01110,1111111,101,1001\}.
$$
In this case $\cU_{\cB}$ has minimum $3$ and maximum $\infty$ repetitions of symbol $0$, and minimum $4$ and maximum $6$ repetitions of symbol $1$. Note that this subshift can be expressed as a 1-step subshift on blocks of $7$ consecutive symbols, which in its most general form will have $2^7$ states and at most two arrows from each state. One can represent $\cU_{\cB}$ using a smaller number of states using multiple representations of the same state, this is shown in Figure~\ref{fig:minmaxexample}.

If the subshift $\cV$ supports an ergodic shift-invariant measure then one can apply tools from ergodic theory \cite{lind2021introduction,adler1992symbolic} or random dynamical systems \cite{arnold1995random}. In particular, if sequences are chosen with respect to an ergodic shift-invariant probability measure $\mu$ on $\cU$ then one can use this to ignore annoying sequences in $\cU$ as long as they lie on some set the can be shown to be zero-measure. 

For example, given any set of non-zero probabilities $P$ assigned to symbols $U$, the Bernoulli product measure $\mu=P^{\bZ}$ assumes the probability of a symbol is given, independent of location, by the same distribution $P$. This can be shown to be an ergodic measure that is zero for any set of sequences where the frequency of each symbol does not match the probability $P$.

\begin{figure}
	\centering
	\includegraphics[width=8cm]{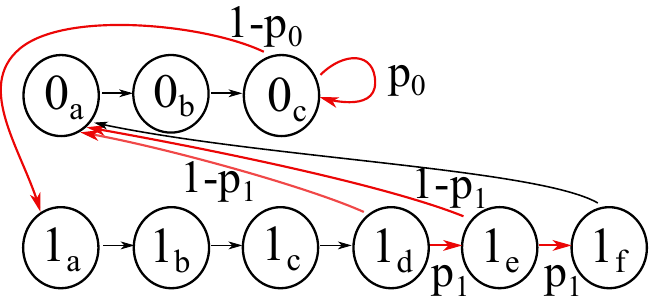}
	\caption{A Markov process for the topological subshift in Figure~\ref{fig:minmaxexample}. For any choice of parameters $0<p_0<1$ and $0<p_1<1$ the transition probabilities shown generate a Markov process with min-max repetitions $m_0^-=3$, $m_0^+=\infty$, $m_1^-=4$ and $m_1^+=6$. For example, after in initial 3 repeats, symbol $0$ will be repeated with probability $p_0$ for each further repeat.}
	\label{fig:minmaxexample_markov}
\end{figure}

Note that there are many ergodic invariant measures corresponding to the topological subshift shown in Figure~\ref{fig:minmaxexample}. By allocating transition probabilities one can define a Markov process that explores this subshift according to this measure. A useful example of this is the family of ergodic subshifts on two symbols 
$$
\cU(m_0^{\pm},m_1^{\pm},p_0,p_1)
$$
with minimum and maximum repeats $m_i^{\pm}$ and repeat probability $0<p_i<1$, for each additional repeat of $i$ after the minimum number. Figure~\ref{fig:minmaxexample_markov} illustrates an example of a family of measures that correspond to the topological subshift in Figure~\ref{fig:minmaxexample} and depends on the parameters $p_0$ and $p_1$.

\subsection{Minimum echo for multiple attractor maps}
\label{subsec:minecho}

We now discuss some testable (but restrictive) assumptions on the maps (\ref{eq:switching_dynamics}) that can be used to give bounds on echo index. Suppose $X \subset \bR^{n}$ is a compact $n$-dimensional manifold, write $\ell(\cdot)$ to denote Lebesgue measure on this and consider the finite set of $M$ maps $f_i: X\rightarrow X$ for $i=0, \ldots,M-1$. 

\begin{assum}
\label{assum:uasps}
Suppose that for each $i=0,\ldots, M-1$ we have:
\begin{itemize}
    \item[(i)] $f_i$ is a continuously differentiable map that is almost everywhere a local diffeomorphism.
    \item[(ii)] $f_i$ has a finite number of hyperbolic stable fixed points $x_i^{0}, \ldots, x_i^{L(i)-1}$.
    \item[(iii)] The basins of attraction $\cB_i^{j}$ of $x_i^j$ exhaust the measure of state space in measure, i.e. 
    $$
    \ell( X \setminus \bigcup_j \cB_i^{j} ) = 0.
    $$
    \item[(iv)] There is a $P(i,j,k)\in\{1,\ldots,L(k)\}$ such that 
    $$
    	%\label{eq:basin}
    	x_i^j\in \cB_k^{P(i,j,k)}.
    $$
    for all $k=1,\ldots,M-1$ and $j=0,\ldots,L(i)-1$.
\end{itemize}
\end{assum}

\begin{remark}
	Condition (i) implies in particular that the pre-image of any set of zero Lebesgue measure under $f_i$ also has zero measure. Condition (ii) implies that near each $x_i^j$ some iterate of $f_i$ is locally a contraction in some neighbourhood; we characterise this in Lemma~\ref{lem:iteratecontraction}. Condition (iii) means that $\{\cB_i^j\}_{j=1}^{L(i)}$ is a full measure partition for each $i$. Condition (iv) means that $x_i^j$ lies within the basin of an attractor for $f_k$ but moreover it is a non-degeneracy assumption that means that no attractor for $f_i$ is on the basin boundary for $f_k$. 
\end{remark}

We characterise the local contraction property more precisely:

\begin{lem}
	\label{lem:iteratecontraction}
	Suppose that Assumption~\ref{assum:uasps} is satisfied. Then for each $i,j$ and any choice of $0<\rho<1$ there is a neighbourhood $N_{i}^j$ of $x_i^j$ and $n_i^j\geq 1$ such that if
	$$
	F:=(f_i|_{N_i^j})^{k}
	$$
	for any $k\geq n_i^j$
	then $F:N_i^j\rightarrow N_i^j$ has a unique fixed point at $x_i^j$ and moreover $F$ contracts by $\rho$
	$$
	\|F(x)-F(y)\|<\rho\|x-y\|.
	$$
\end{lem}

\proof The assumption of a hyperbolic attraction fixed point $x_i^j$ implies linear stability, and hence that in any neighbourhood contained in $\cB_i^j$ some iterate of $f_i$ will be a contraction.
\qed

Now consider a specific input sequence $\bv=\{v[k]\}_{k\in\bZ}$. Given some choice of $A[0]\in\{1,\ldots,N_{v[0]}\}$ and applying Assumption~\ref{assum:uasps}(iv) we have
$$
x_{v[0]}^{A[0]}\in \cB_{v[1]}^{P(v[0],A[0],v[1])}.
$$
Hence, associated with a sequence $\bv$ and initial choice of attractor $x_{v[0]}^{A[0]}$ there will be a unique sequence
$$
\left\{x_{v[k]}^{A[k]}\right \}_{k\geq 0}
$$
such that  $x_{v[k]}^{A[k]}$ is an attractor for $f_{v[k]}$ contained in the basin of $x_{v[k+1]}^{A[k+1]}$ for the map $f_{v[k+1]}$, namely
\begin{equation}
	\label{eq:attractorseq}
	A[k]=P(v[k-1],A[k-1],v[k]).
\end{equation}
We call such a sequence $A[k]$ a {\em forward attractor sequence} for $\bv$ starting at $A[0]$. Since $A[k]$ is determined by (\ref{eq:attractorseq}) there will be only finitely many forward attractor sequences and the number of these is bounded above by $\min_{i} L(i)$. An {\em entire attractor sequence} is a sequence $A[k]$ satisfying (\ref{eq:attractorseq}) that is associated with a bi-infinite $\bv$; this depends not only on the $v[k]$ for $k\geq 0$ but also on $k<0$.

\begin{lem}
\label{lem:entire}
Suppose that Assumption~\ref{assum:uasps} holds for the system (\ref{eq:switching_dynamics}). Then there is a $m_{\min}$ such that for any $\bv$ with $m_i^-\geq m_{\min}$ for all $i$ and any entire attractor sequence $A[k]$ for this $\bv$, there is a pullback attracting entire solution $x^{A[k]}_{v[k]}$ such that 
$x^{A[k]}_{v[k]}\in \cB_{v[k+1]}^{A[k+1]}$.
\end{lem}
	
\proof
Choose any $0<\rho<1$; by applying Lemma~\ref{lem:iteratecontraction} for all choices of $i,j$ we can find an $\epsilon>0$ and an $m_{\min}$ such that $m_{\min}>n_i^j$ and $B_{\epsilon}(x_i^j) \subset N_i^j$. Pick any attractor sequence $A[k]$ for a $\bv$ that satisfies $m_i^-\geq m_{\min}$ and define
$$
x^{[A]}[k,n]:= \Phi_{n,\sigma^{k-n}(\bv)} B_{\epsilon}(x_{v[k-n]}^{A[k-n]}).
$$
This is a nested sequence in increasing $n$ for fixed $k$ and there is contraction by $\rho$ over blocks of the same symbol. Hence the set is non-empty and has diameter that shrinks to zero as $n\rightarrow \infty$. This means that 
$$
x^{[A]}[k]:= \bigcap_{n<k} x^{[A]}[k,n]
$$
consists of a single entire solution that pullback attracts an $\epsilon$-neighbourhood of itself.
\qed

A consequence of this is that, for long enough minimum block-lengths we can get a lower bound for echo index from the number of distinct entire attractor sequences.

\begin{thm}
	\label{thm:suff}
Suppose that Assumption~\ref{assum:uasps} holds for (\ref{eq:switching_dynamics}) and choose $m_{\min}$ and $\bv$ such that the conclusion of Lemma~\ref{lem:entire} holds. Suppose that there are $E$ distinct entire attractor sequences for $\bv$. Then the echo index is at least $E$ for this sequence.
\end{thm}

\proof
By Lemma~\ref{lem:entire} each entire attractor sequence $A[k]$ has a pullback attractor $x^{[A]}[k]$. These are distinct as long as the entire attractor sequences are distinct.
\qed

Under additional assumptions, one can ensure echo index, in particular for the following case where we assume there is a uniform bound on how long it takes points to enter a given neighbourhood of a single attractor.

\begin{thm}
\label{thm:echoone}
	Suppose that Assumption~\ref{assum:uasps} holds for (\ref{eq:switching_dynamics}) and in addition assume that there is an $i$ such that $f_i$ has a single attracting fixed point $x_i^0$ such that $h(f_i^m(X),x_i^0)\rightarrow 0$ as $m\rightarrow \infty$. Then there is a single attractor sequence and an $m_{\min}$ such that for every $\bv$ with $m_i^-\geq m_{\min}$ the system has echo index one.
\end{thm}

\proof
In this case note that whenever $v[k]=i$ we have $A[k]=0$. Hence there is only one entire attractor sequence. Moreover, given any $\epsilon$ there is an $m$ such that all points must enter the neighbourhood $B_{\epsilon}(x_i^0)$ after at most $m$ iterates and hence the entire solution pullback attracts all points in $X$.
\qed

\subsection{Obstructions to conditions for echo consistency}
	
Theorem~\ref{thm:suff} gives sufficient conditions to guarantee a minimum echo index for the forced system in terms of $E$ the number of distinct attractor sequences for the system with input $\bv$. Conversely, Theorem~\ref{thm:echoone} shows under stronger conditions that the echo index is precisely one.

One might naively expect that an even stronger result than Theorem~\ref{thm:suff} may follow, namely that the echo index is in fact $E$ for cases where $E\geq 2$. This is not the case because under composition new attractors may appear near a basin boundary. 

\begin{figure}
	\centering
\includegraphics[width=12cm]{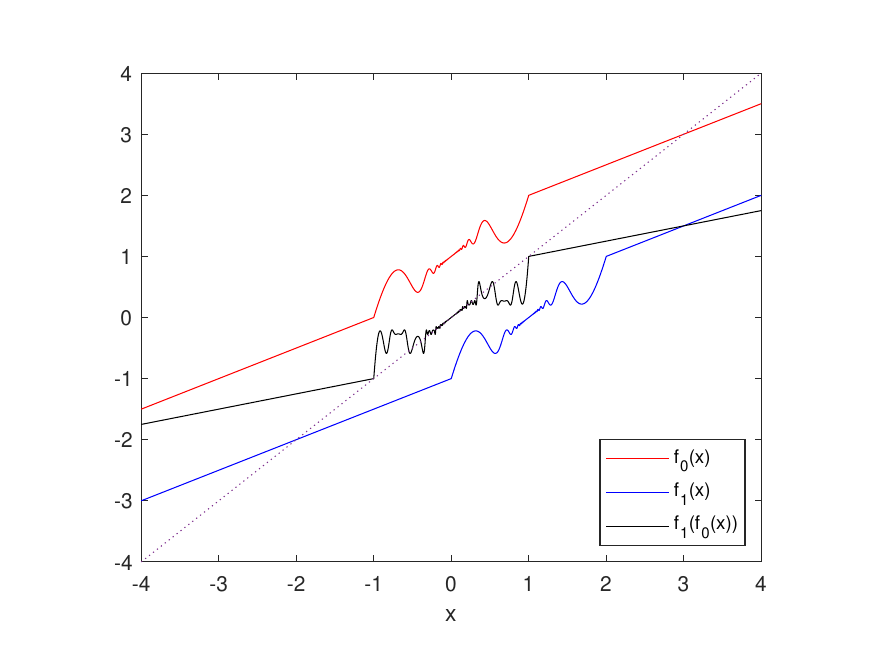}
\caption{Maps $f_0$ and $f_1$ (\ref{eq:f0f1}) such that each individual map has a unique globally attracting linearly stable fixed point. However the composition of $f_1\circ f_0$ has infinitely many stable fixed points.}
\label{fig:diabolic}
\end{figure}

As an example, define
$$
F(x)=x+x^2\sin\frac{\pi}{x},~~F(0)=0
$$
and $f(x)$ is the continuous function such that $f(x)=F(x)$ for $|F(x)|<1$ and $f(x)$ is a linear function with constant slope $0.5$ for $|F(x)|>1$. Now consider the maps
\begin{equation}
	\label{eq:f0f1}
f_0(x) = f(x)+1,~~~
f_1(x) = f(x-1).
\end{equation}
One can verify that $f_1(x)$ and $f_2(x)$ each has a unique attracting fixed point, namely
$x_0^0=3$, $x_1^0=2$. However, $f_0\circ f_1=f^2$ has infinitely many attracting fixed points separated by repelling fixed points; it also has neutrally stable and linearly unstable fixed points; see Figure~\ref{fig:diabolic}. In summary, this system satisfies Assumption~\ref{assum:uasps} and for any input there is only one attractor sequence $A[k]=0$ for all $k$, but for an input $\bv$ that alternates between $0$ and $1$ there are infinitely many attracting entire solutions.

Hence for this input there is infinite echo index even though each. Nonetheless, applying Theorem~\ref{thm:echoone} there is a minimum block-length such that there is only one attracting entire solution. In this case, numerical simulations suggest this minimum is $2$, which would suggest that the system has echo index one for all input sequences except the periodic sequence that is an infinite repeat of $01$.

\section{Echo index dependence for RNNs}
\label{sec:rnns}

One of the main motivations for this study is the need to understand how many responses there are for RNNs such as (a) trained discrete-time RNNs of the form \cite{tallec2018can,bianchi2017recurrent}
\begin{align}
	\label{eq:rnn}
	x[k+1] & = \phi( W_r x[k] + W_i u[k+1] + W_{f} z[k] )
\end{align}
driven by inputs $u[k]$ or (b) trained ESNs of leaky-integrator neurons:
\begin{equation}
	\label{eq:leaky_rnn}
	x[k+1] = G(x[k],u[k],z[k]),~~G(x,u,z)=(1-\alpha) x + \alpha \phi( W_r x + W_i u + W_{f} z ),
\end{equation}
where $\alpha \in (0,1) $ quantifies the leakage, $\phi$ is a nonlinear function; in both case the input sequence is given by $u[k]$ and $z[k]$ represents the output feedback
\begin{equation}
		\label{eq:output}
	z[k] = \psi(x[k]).
\end{equation}
The nonlinear function $\phi$ is called an \emph{activation function}, we assume it is bounded, monotonically increasing, and differentiable, e.g. a scaling of $\tanh$. By contrast, $\psi$ is usually the identity function or a softmax function.

\subsection{Switching dynamics for RNNs}

We consider an ESN with leaky-integrator neurons \cite{jaeger2007optimization} and no output feedback\footnote{As discussed in \cite{ceni2020interpreting} whenever the readout is linear the feedback term can be formally incorporated in the reservoir term. Therefore the absence of the feedback term in an ESN state-update rule can represent an ESN state-update rule after the training session where $W_r$ represents the ``effective reservoir'' after training.} for implementing the input-driven state update rule \eqref{eq:switching_dynamics}.  We consider a map of the form (\ref{eq:leaky_rnn}) but with
\begin{equation}
\label{eq:map}
G_{\alpha}(x, u):= (1-\alpha)x + \alpha \tanh(W_r x + W_{in} u).
\end{equation}
A finite set of input values $U = \{ u_0, u_1, \ldots, u_{M-1} \}$ defines a number $M$ of autonomous maps $f_i: X \rightarrow X $, where $X=[-1,1]^r$ and $r$ the dimension of the internal state of the RNN, i.e. the number of neurons.
In this RNN's framework, we can see that for small enough input values, the echo index of the nonautonomous switching system is the number of attractors of the input-free ESN.
On the other hand, Theorem~\ref{thm:echoone} has an interpretation in terms of large amplitude forcing for RNN-like systems. In fact, we know that large amplitude inputs drives the system in the saturating tails regime of $\tanh$ characterised by a single attracting fixed point \cite{ceni2020interpreting}. Thus Theorem~\ref{thm:echoone} implies that on forcing an RNN with a large amplitude input for long enough, the resulting nonautonomous RNN switching dynamics are characterised by echo index one.

\subsection{An example of input-driven RNN with multiple attractors}
\label{sec:example}

We provide here an example with a two dimensional ESN to better illustrate the concepts. 
We choose,
\begin{equation}
    W_r = \begin{bmatrix}
    \frac{1}{2}  &  0 \\
    0     & \frac{7}{4}
\end{bmatrix},
\quad W_{in} = I_2,
\end{equation}
with $I_2$ the identity matrix.  We consider input sequences $\cU = \{ u_0, u_1  \}^{\bZ}$ where 
\begin{equation}
\label{eq:inputs}
u_0:=  \begin{pmatrix}     
    \frac{1}{4}  \\
    \frac{1}{20}
\end{pmatrix}~\mbox{ and }~u_1:=  \begin{pmatrix}  
    -\frac{1}{4}  \\
    -\frac{1}{2}
\end{pmatrix}.
\end{equation}
What follows can be observed for any value of $ \alpha \in (0,1] $. 
We chose $\alpha=\dfrac{1}{4}$ again because a small leak rate highlights the transient dynamics. 

\begin{figure}[ht!]
	\centering
 \includegraphics[width=16cm]{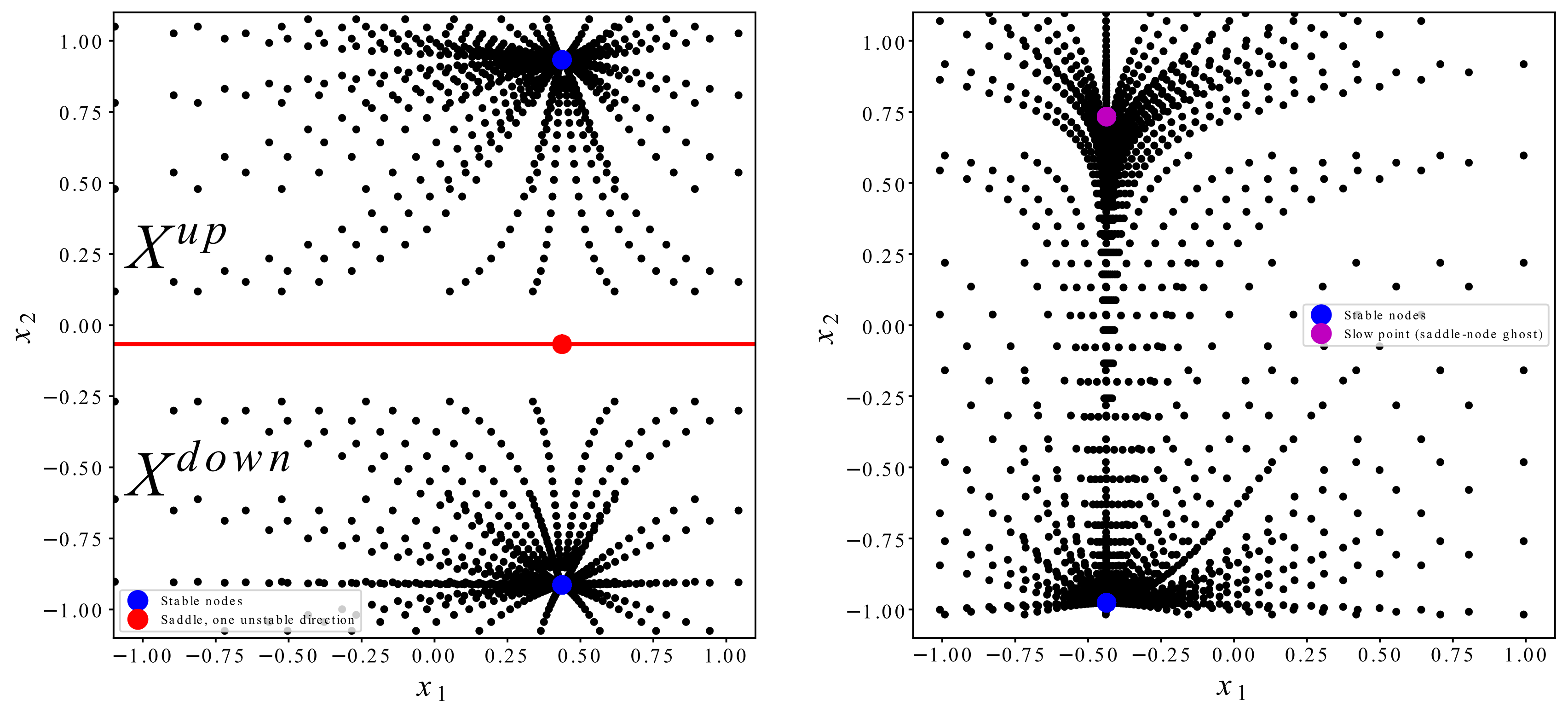}

 (a) \hspace{7cm} (b)
	\caption{
		(a) Phase portrait of the map $f_0$. Red line represents the stable manifold of the saddle. 
		Some initial conditions have been evolved and plotted (as black points) in order to visualise the vector field.
		(b) phase space portrait of the map $f_1$. Note that the purple point is not a fixed point but a \emph{slow point} \cite{sussillo2013opening,ceni2020interpreting}. In fact, the map $f_1$ is close to a saddle-node bifurcation which occurs nearby the position of such a slow point.
	}
	\label{fig:autonomous_maps}
\end{figure}

The nonautonomous dynamics driven by some input sequence $\bu \in \cU$ consists of sequence of applications of the two maps: 
$$
f_0(x):= G(x,u_0),~~f_1(x):= G(x,u_1).
$$
One can verify that the autonomous system $x[k+1] = f_0(x[k])$ has two asymptotically stable points with a saddle between them along the vertical line of $x_1 \approx 0.45 $ (see Figure \ref{fig:autonomous_maps}) while the autonomous system $x[k+1] = f_1(x[k])$ has only one (asymptotically stable) fixed point lying in the quadrant where both variables are negatives. Note that \cite[Theorem 4.1]{ceni2020echo} can be applied in this example to prove the existence of a local point attractor lying in a strip of negative values of the $x_2$ variable. Nevertheless, it is not straightforward to prove the existence of additional local point attractors.

The stable manifold of the saddle of the autonomous map $f_0$ is a horizontal line dividing the phase space in two sets.
Let us denote with $x^*$ the upper stable node of $f_0$ lying on the quadrant where both variables are positive. Thanks to \cite[Proposition B.1]{ceni2020echo}, we can consider the phase space to be $X=[-1,1]^2$.
Let us call $ X^{up}$ the upper half where $x^*$ lies (including the stable manifold line) and $ X^{down} $ the remaining part.
On the other hand, the global attractor of the autonomous map $f_1$ consists of only an asymptotically stable fixed point lying in $ X^{down}$. Lemma~\ref{lem:iteratecontraction} can be applied to show there exists a (minimum) positive integer $m_{\min}$ for which $ f_1^{m_{\min}}( X^{up} ) $ is mapped into the interior of $X^{down}$. For the particular choice of parameters it turns out that $m_{\min}\approx 30$.

\subsection{Bifurcations of echo index}
\label{sec:bif-2esp-1esp}

\begin{figure}
\includegraphics[width=8cm]{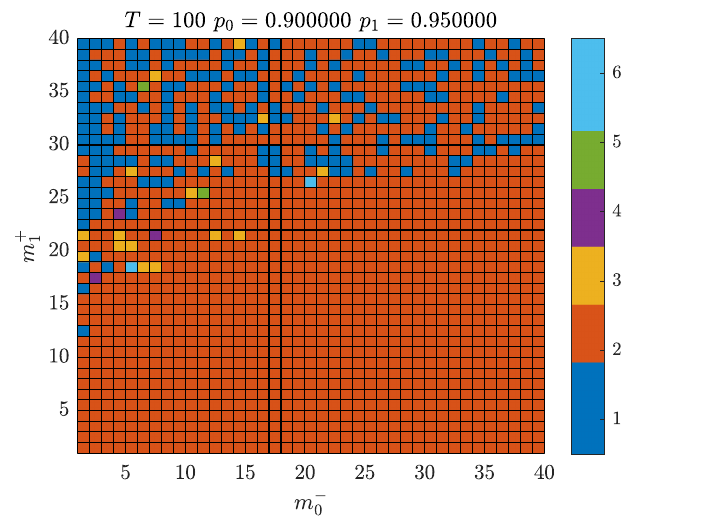}
~\includegraphics[width=8cm]{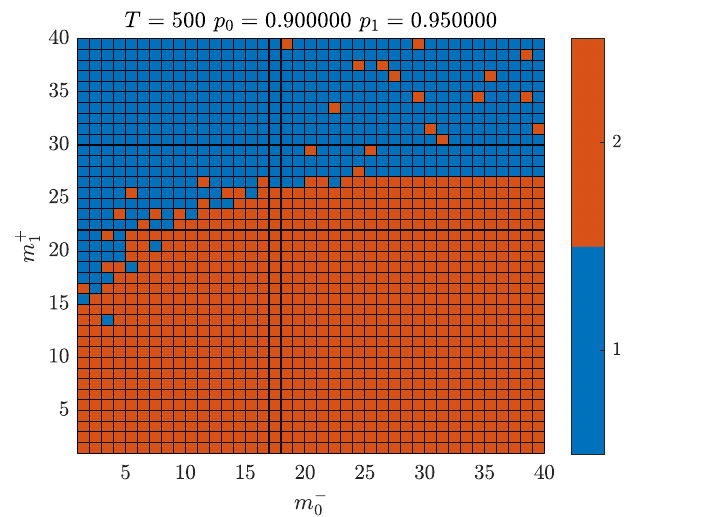}

(a)\hspace{8cm}(b)

~

\includegraphics[width=8cm]{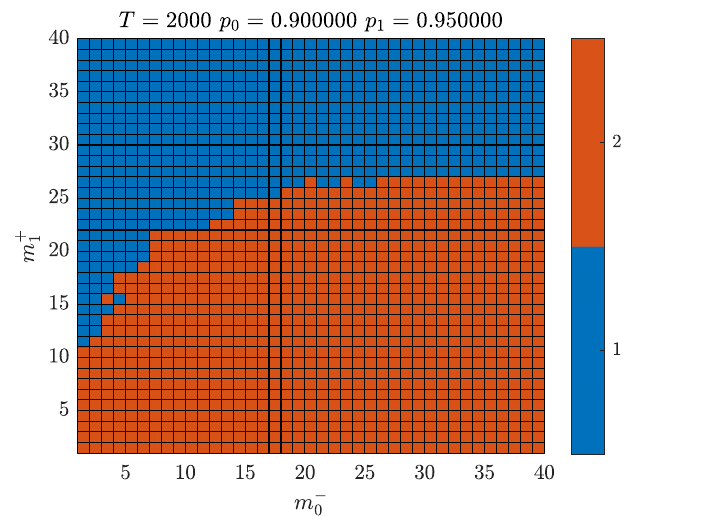}
~\includegraphics[width=8cm]{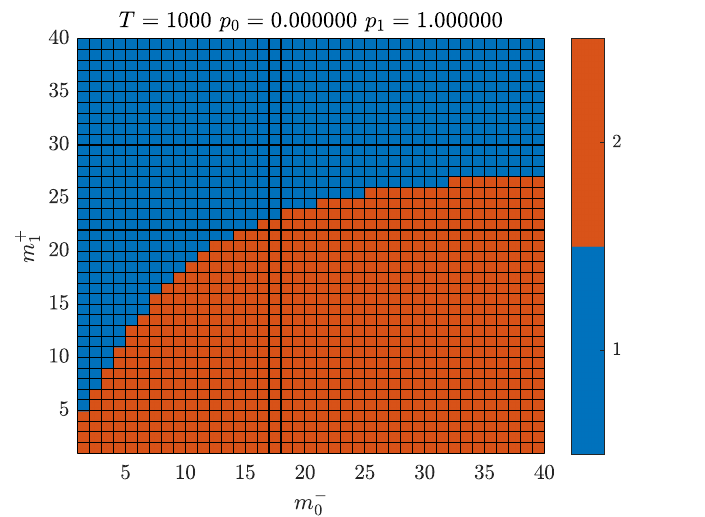}

(c)\hspace{8cm}(d)
\caption{Estimates of echo index for the system (\ref{eq:map}) with parameters as in the text for a range of random input sequences.  Random sequences of length $T$ with varying $m_0^-$ and $m_1^+$ are chosen, fixing  $m_0^+=40$ and $m_1^-=1$. Note that the echo index is apparently 2 or more for small enough $m_1^+$, or for small $T$. (a-c) show estimates of the echo index for randomly generate inputs with the given min-max block-lengths and probabilities of repetition $p_0$ and $p_1$ as in Figure~\ref{fig:minmaxexample_markov}. Examining longer timeseries when estimating the echo index leads to a clear boundary between regions of different echo index. (d) shows a special case where there are only periodic inputs. We conjecture that the cases (a-c) will limit to (d) for arbitrarily large $T$. 
}
\label{fig:random_inputs}
\end{figure}

In this section we perform numerical simulations to compute the echo index of the ESN's example of section~\ref{sec:example}.
For each choice of parameters determining the sequence $\bv$, we choose $50$ uniformly distributed initial conditions and iterate $T$ steps before using a clustering algorithm to numerically estimate the number of clusters in the final state - this is an estimate of the echo index. Random sequences of length $T$ with varying $m_0^-$ and $m_1^+$ are chosen, and we fix $m_0^+=40$ and $m_1^-=1$. Figure~\ref{fig:random_inputs} highlights that for $m_0^+$ large enough there will be echo index $2$ for small enough $m_1^+$. Panel (a) shows $T=100$, $p_0=0.9$ and $p_1=0.95$. Note that the short length of the timeseries is not enough to collapse the initial conditions down to one of two values. Panel (b) is as for (a) but with $T=1000$; in this case we have index one or two, though some cases where $m_1^+$ is large have not been sampled for long enough giving some spurious echo index $2$. Panel (c) is as for (a) but with $T=2000$; we find what is presumably a clear boundary emerging between different values of the echo index. Finally (d) shows the deterministic case where $p_0=0$ and $p_1=1$ which corresponds to the periodic orbit where there are $m_0^-$ repeats of $0$ and $m_1^+$ repeats of $1$ - in this case there a clear boundary between regions with echo index one and two; this presumably corresponds to a bifurcation of the periodically forced map where a second attractor appears.  Note that for $m_1^+>30$ and large enough $T$, we always find echo index one as suggested by the discussion at the end of Section~\ref{sec:example}, but for $m_0^-$ small it is possible to find echo index one for some values of $m_1^+<30$.

\section{Conclusions} 
\label{sec:conclusions}

In this paper we have gone beyond work in \cite{ceni2020echo} to highlight specific ways in which the echo index varies with input signal (and hence whether the echo state property holds for a given input). We present this for an iterated function system on a compact space with an example application to an echo state network.

We have so far only considered min-max block-length and repetition probabilities on determining bounds on echo index but the actual value may depend on much more subtle properties of words appearing in the input and properties of the individual maps. It remains a challenge to better understand this relationship between input set, system properties and echo index, even if we restrict only to functions where the only attractors are fixed points. Clearly, responses for cases where the autonomous dynamics of the maps include more complex attractors (such as chaotic, quasiperiodic or period) will be more challenging, not least because the a simple generalization of Assumption~\ref{assum:uasps}(iv) is unreasonable - a single attractor of one map will can stably straddle several basins of attractions for attractors of another map.

One interpretation of our results is that minimum block-lengths are a proxy for the input rate to the system - a long enough minimum block-length corresponds to a slow rate of input. Our results Theorems~\ref{thm:suff} and \ref{thm:echoone} imply that the expected behaviour can be characterised by attractors and basins of the individual maps for slow enough rates of input. For shorter min block-lengths the picture becomes more complex and the transient or nonautonomous behaviour of each map features more strongly in the response to input - there will be an analogy to rate-induced critical transitions \cite{ashwin2012tipping} in the response for such cases.

A future direction that seems worthy of study is the role of transients in determining the echo index. The computations in Figure~\ref{fig:random_inputs} show that even for quite long (but finite) computations the number of responses may apparently exceed the echo index which is attained asymptotically. A thorough understanding of responses will be needed to explain such transient behaviour of the echo index. This clearly depends not only on transients in the map dynamics but also on waiting times to see certain words within the inputs. 
% We do not deeply consider applications other than to highlight that echo index greater than one is not necessarily an indication that the trained RNN is malfunctioning for this input - it could also be used to design multifunction RNNs \cite{ceni2021phd,multifunction2022recurrent} that will remember their task.

\subsection*{Acknowledgements}

We thank EPSRC for support via EP/W52265X/1. We thank Lorenzo Livi, Muhammed Fadera and Claire Postlethwaite for very informative discussions in relation to this work.

\subsection*{Data Access}

The Matlab code for Figure~\ref{fig:random_inputs} is available from {\tt https://github.com/peterashwin/ashwin-ceni-2023}.

~~~

%\section*{References}

%\bibliographystyle{elsarticle-num}
\bibliographystyle{plainnat}
\bibliography{nesp-refs-dec22}

\end{document}